\documentclass[12pt,reqno]{amsart}

\usepackage[T1]{fontenc}
\pdfoutput=1
\usepackage{mathptmx,microtype,graphicx,dsfont}
\usepackage[hidelinks]{hyperref}
\usepackage{amsthm,amsmath,amssymb}
\usepackage{enumitem}
\usepackage[nameinlink,capitalise]{cleveref}
\usepackage[font=footnotesize,margin=2cm]{caption}
\usepackage{cite}

\usepackage[foot]{amsaddr}

\crefname{theorem}{Theorem}{Theorems}
\crefname{thm}{Theorem}{Theorems}
\crefname{mainthm}{Theorem}{Theorems}
\crefname{lemma}{Lemma}{Lemmas}
\crefname{lem}{Lemma}{Lemmas}
\crefname{remark}{Remark}{Remarks}
\crefname{claim}{Claim}{Claims}
\crefname{prop}{Proposition}{Propositions}
\crefname{proposition}{Proposition}{Propositions}
\crefname{defn}{Definition}{Definitions}
\crefname{corollary}{Corollary}{Corollaries}
\crefname{conjecture}{Conjecture}{Conjectures}
\crefname{question}{Question}{Questions}
\crefname{section}{Section}{Sections}
\crefname{figure}{Figure}{Figures}
\crefname{table}{Table}{Tables}
\usepackage{xcolor}

\theoremstyle{plain}

\newtheorem{thm}{Theorem}
\newtheorem*{thm*}{Theorem}
\newtheorem{lemma}[thm]{Lemma}

\theoremstyle{definition}

\theoremstyle{remark}

\numberwithin{equation}{section}

\newcommand{\cI}{{\mathcal I}}

\newcommand\Z{\mathbb{Z}}
\newcommand\R{\mathbb{R}}

\newcommand\E{{\mathbb E}}

\renewcommand{\le}{\leqslant}
\renewcommand{\ge}{\geqslant}

\author[M. Bassan]{Michal Bassan}
\address{University of Oxford, 
Department of Statistics
and Keble College}
\email{michal.bassan@keble.ox.ac.uk}

\author[S. Donderwinkel]{Serte Donderwinkel}
\address{University of Groningen, 
Bernoulli Institute for Mathematics, 
Computer Science and AI, 
and CogniGron 
(Groningen Cognitive Systems and Materials Center)}
\email{s.a.donderwinkel@rug.nl}

\author[B. Kolesnik]{Brett Kolesnik}
\address{University of Warwick, 
Department of Statistics}
\email{brett.kolesnik@warwick.ac.uk}


\begin{document}

\title[Score sequences, EGZ numbers, and the LK method]
{Tournament score sequences,  
Erd\H{o}s--Ginzburg--Ziv numbers,
and the L\'evy--Khintchine method}

\begin{abstract} 
We give a short proof 
of a recent result of 
Claesson, Dukes, Frankl\'in and Stef\'ansson, 
connecting the number $S_n$ of  
score sequences 
and the Erd\H{o}s--Ginzburg--Ziv numbers $N_n$ 
from additive number theory. 
Our proof utilizes the lattice path 
representation of 
score sequences
by 
Erd\H{o}s and Moser, and remarks
by Kleitman added to an article of Moser
regarding cyclic shifts of such paths.
The  connection 
between $S_n$ and $N_n$
is an instance of 
the L\'evy--Khintchine 
formula from probability
theory. 
We highlight the utility of such formulas,  
by giving a short proof 
of Moser's conjecture that 
$S_n\sim C4^n/n^{5/2}$, 
where $C$ is described in terms of $N_n$. 
\end{abstract}

\maketitle

\section{Introduction}\label{S_intro}

A {\it tournament} 
is an orientation of the complete graph 
$K_n$. We think of vertices as players and edges as games, 
with each edge directed towards the winner. 
The {\it score sequence} 
lists the total number of wins by the players
in non-decreasing order.
Landau \cite{Lan53} showed that 
$s_1\le\cdots\le s_n$ in $\Z^n$
is a score sequence if 
$\sum_{i=1}^n s_i={n\choose2}$, with all partial sums
$\sum_{i=1}^k s_i\ge{k\choose2}$. 
The conditions are necessary, 
since any $k$ teams play ${k\choose2}$
games amongst themselves.

In this work, we give a short proof 
of Moser's \cite{Mos71} conjecture 
that the number $S_n$ of score sequences
satisfies $S_n\sim C 4^n/n^{5/2}$. 
In doing so, we will highlight a novel 
probabilistic method of asymptotic enumeration, 
which we 
expect to find more combinatorial applications.

This method 
is founded on the study of 
analytic transformations $\Phi(\mu)$ of probability 
distributions $\mu$ by  
Chover, Ney and Wainger \cite{CNW73}. 
Classical renewal theory 
corresponds to the special case that $\Phi(z)=(1-z)^{-1}$, 
in which case $\Phi(\mu)=\sum_{m=0}^\infty \mu^{*m}$ 
is a sum over the convolutions of $\mu$. 
Using results in \cite{CNW73}, Hawkes and Jenkins \cite{HJ78}
(cf.\ Embrechts and Hawkes \cite{EH82}) 
obtained conditions under which the asymptotics
of a sequence $A_n$ and those of a certain transform $A_n^*$ 
are related
as $A_n/A_n^*\sim C/n$, for some constant $C$. 
The specific case $C=1$ was analyzed earlier
by Wright \cite{Wri67a,Wri67b,Wri67c}.

The sequences $A_n$ and $A_n^*$ 
and the constant $C$ 
can be described in terms of the 
{\it L\'evy--Khintchine formula} from probability theory;  
see \cref{S_LKmethod}. 
The power of this method, based on our recent 
experience \cite{Kol23,DK24b,BDK24}, 
is that $A_n^*$ can be much simpler than $A_n$.

Claesson, Dukes, Frankl\'in and Stef\'ansson \cite{CDFS23} 
recently proved that 
\begin{equation}\label{E_CDFS}
nS_n 
= \sum_{k=1}^n N_k S_{n-k},
\quad\quad
n\ge1,
\end{equation}
where $N_n$ is the number of 
subsets of $\{1,\dots, 2n-1\}$ of size $n$ 
whose elements sum to $0$ mod $n$. 
We call $N_n$ the {\it Erd\H{o}s--Ginzburg--Ziv numbers,} 
with reference to their result \cite{EGZ61}
that {\it any} set of $2n-1$ 
integers has such a subset. 

As discussed in \cref{S_LKmethod}, 
\eqref{E_CDFS}
implies that $S_n^*=N_n$. 

In the early 1900s, 
von Sterneck \cite{Bac02} 
(cf.\ \cite{Ram44,Ale08,Che19}) 
showed that
\begin{equation}\label{E_Nn}
N_n = 
\sum_{d|n}\frac{(-1)^{n+d}}{2n}{2d\choose d}\phi(n/d), 
\end{equation}
where $\phi$ is Euler's totient function. 

In \cite{Kol23}, the third author observed that,  
by combining \eqref{E_CDFS} and \eqref{E_Nn}
with the limit theory in \cite{HJ78}, 
it follows that that $S_n\sim C 4^n/n^{5/2}$, 
as conjectured by Moser \cite{Mos71}. 

\begin{thm}
\label{T_main}
As $n\to\infty$, 
we have that 
\begin{equation}\label{E_asySn}
\frac{n^{5/2}}{4^n}S_n
\to \frac{1}{2\sqrt{\pi}}
\exp \left(
\sum_{k=1}^\infty\frac{N_k}{k4^k}
\right). 
\end{equation}
\end{thm}

\subsection{Purpose}
In this work, we give 
a simple proof of 
\eqref{E_CDFS} using: 
\begin{enumerate}
\item the lattice path representation of 
score sequences, first observed by Erd\H{o}s and Moser, and 
\item Kleitman's brief remarks, added
to the end of 
Moser's article \cite{Mos71}, 
regarding cyclic
shifts of score sequences. 
\end{enumerate}
As a result, we obtain a short
proof and deeper explanation
for \cref{T_main}, based also 
on the probabilistic point
of view discussed in \cref{S_LKmethod}
below.

The proof of \eqref{E_CDFS} in \cite{CDFS23}
is more involved, 
as it takes place entirely at the level of 
modular arithmetic and generating functions. 

As discussed in 
Moon \cite[Theorem 33]{Moo68}, 
the relationship between score sequences
$s_1\le\cdots\le s_n$ 
of length $n$
and up/right lattice paths 
from $(0,0)$ to $(n,n)$ 
goes back to Erd\H{o}s and Moser
in the 1960s. Informally, consider 
the bar graph of the sequence, where 
the $i$th bar has height $s_i$. 
Rotating such a lattice path gives
an up/down bridge of length $2n$.

We prove \eqref{E_CDFS} using the
renewal structure of $S_n$   
(see, e.g., Moon \cite[\S2--3]{Moo68}).
We decompose bridges 
associated with score sequences 
into irreducible parts, 
and argue that cyclic shifts are related
to the Erd\H{o}s--Ginzburg--Ziv numbers $N_n$.
In \cref{S_areas}, we introduce 
the {\it diamond area} $a(B)$ of a 
bridge $B$, which 
reveals 
the geometric connection between $S_n$
and $N_n$.

Kleitman observed that $S_n$ 
can be bounded by 
considering 
cyclic shifts of bridges $B$ with $a(B)=0$,
as noted in \cref{S_Kleitman}.
Building on this, in \cref{S_proof,S_OurCase}, 
we show that cyclic shifts of bridges $B$ with $a(B)\equiv 0$
mod $n$ are counted by $N_n$ and lead to the precise asymptotics of $S_n$.

\subsection{Combinatorial geometry}
\label{E_perm}
\begin{figure}[h!]
\centering
\includegraphics[scale=0.65]{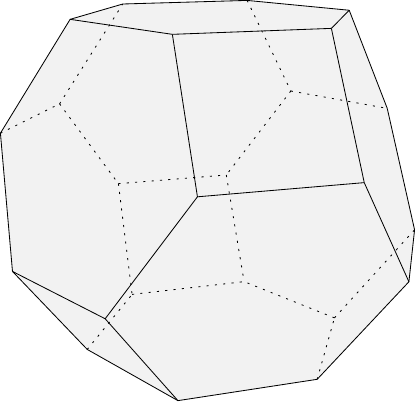}
\caption{The permutahedron 
$\Pi_3\subset\R^4$ (projected into $\R^3$) 
is the convex hull of 0123 and its permutations. 
Its non-decreasing lattice
points 0123, 0222, 1113 and 1122
are the 
$S_4=4$
score sequences. 
}
\label{F_per}
\end{figure}

The permutahedron $\Pi_{n-1}$ is 
a classical object in discrete 
geometry, obtained as 
the convex hull
of the score sequence $(0,1,\ldots,n-1)$ and its permutations; see
\cref{F_per}. 
Score sequences correspond to its
non-decreasing lattice points;
see, e.g., \cite{KS20}. 
On the other hand, Zaslavsky
observed that the set of {\it all} lattice
points is in bijection with the spanning
forests of $K_n$, as discussed
in Stanley \cite{Sta80}.

Our techniques might be helpful 
with enumerating various classes of 
lattice points in  
the {\it generalized permutahedra} in 
Postnikov \cite{Pos09}
and {\it Coxeter permutahedra} 
in 
Ardila, Castillo, Eur and Postnikov
\cite{ACEP20}.
See \cite{KS20,KS23,KMP23} for connections
between tournaments and these 
more general permutahedra.

\subsection{Acknowledgments}
MB is supported by a Clarendon Fund Scholarship. 
SD acknowledges
the financial support of the CogniGron research center
and the Ubbo Emmius Funds (Univ.\ of Groningen).
BK was supported by a 
Florence Nightingale Bicentennial Fellowship (Oxford Statistics)
and a Senior Demyship (Magdalen College).

\section{The L\'evy--Khintchine method} 
\label{S_LKmethod}

Recurrences of the form \eqref{E_CDFS} 
are related to
the 
{\it L\'evy--Khintchine formula} 
from probability theory (see, e.g., \cite{EH82}). 
This formula characterizes all {\it infinitely divisible} 
random variables $X$ on $\R^d$ (see, e.g., \cite{Ber96}). 
We recall that $X$ is infinitely divisible if for all $n\ge1$, 
there are independent and identically
distributed $X_1,\ldots,X_n$ such that 
$X$ and $X_1+\cdots+X_n$ are equal in distribution. 

We recall that
a positive, summable sequence 
$(1=a_0,a_1,\ldots)$ is proportional 
to an infinitely divisible probability distribution 
$p_n= a_n/\sum_k a_k$ on the integers
$n\ge0$ if and only if 
\begin{equation}\label{E_LK}
\sum_{n=0}^\infty a_nx^n
=\exp\left(
\sum_{k=0}^\infty \frac{a_k^*}{k}x^k
\right)
\end{equation}
for some non-negative sequence 
$(0=a^*_0,a^*_1,\ldots)$.  
See, e.g., \cite{Kat67,Fel68,Ste71,vHa78,HJ78}
for a proof.
Since \eqref{E_LK}
is a special case of the 
L\'evy--Khintchine
formula, we
 call $a^*_n$ the {\it L\'evy--Khintchine
transform} of $a_n$.

Differentiating \eqref{E_LK}
and the comparing coefficients, 
it can be seen that 
\eqref{E_LK} is equivalent to 
the recurrence 
\begin{equation}\label{eq:renewal}
na_n 
= \sum_{k=1}^n a_k^* a_{n-k},
\quad\quad
n\ge1. 
\end{equation}

A positive sequence $\vartheta(n)$ is {\it regularly varying
with index $\gamma$} if, 
for all $x>0$, we have that 
$\vartheta(\lfloor xn\rfloor)/\vartheta(n)\to x^{\gamma}$ 
(see, e.g., Bojanic and Seneta \cite[Corollary 1]{BS73}). 
Hawkes and Jenkins \cite{HJ78} (cf.\ Embrechts and Hawkes \cite{EH82}) 
showed that, 
if $a_n^*$ is regularly 
varying with some index $\gamma<0$, then 
\begin{equation}\label{E_asyAn}
a_n\sim\frac{a_n^*}{n}
\exp\left(
\sum_{k=1}^\infty\frac{a_k^*}{k}
\right). 
\end{equation}

With an eye to applications, we might think of 
$A_n$ as counting 
the size of some class of combinatorial objects. 
Naturally, in this context, if 
$nA_n 
= \sum_{k=1}^n A_k^* A_{n-k}$, 
we call $A_n^*$ 
the L\'evy--Khintchine transform of $A_n$. 
If $A_n$ has exponential growth rate $\alpha$, we 
let $a_n=A_n/\alpha^n$. 
If $a^*_n=A^*_n/\alpha^n$ 
is regularly varying, with some index $\gamma<0$, then 
by \eqref{E_asyAn} we can express the asymptotics of $A_n$ in terms of the sequence 
$(A_1^*,A^*_2,\dots)$.

\section{Transforming renewal sequences}
\label{S_renewal}

{\it Renewal sequences}
are a
special class of sequences $A_n$
that have
L\'evy--Khintchine
transforms $A_n^*$. 
Such sequences arise frequently in combinatorics, 
when counting
structures of length $n$ 
that can be  
decomposed into a series of 
irreducible parts. 
More formally, $A_n$ 
is a renewal sequence if its 
generating function $A(x)=\sum_{n=0}^\infty A_nx^n$ can be expressed as 
\[
A(x)=\frac{1}{1-A^{(1)}(x)},
\]
where $A^{(1)}(x)=\sum_{n=0}^\infty A_n^{(1)}x^n$
is the generating function for the number 
$A_n^{(1)}$ of irreducible structures of length $n$. 
See, e.g., Feller \cite{Fel68}
for details. 

In such cases, 
$A_n^*$ takes a special form,
in terms of cyclic shifts. 

\begin{lemma}
\label{L_renewal}
Suppose that $A_n$ is a renewal sequence. 
Then: 
\begin{enumerate}
\item the L\'evy--Khintchine
transform $A_n^*$ is the number of pairs $(X,m)$, 
where $X$ is a structure of length $n$ and $0\le m <\ell$,
where $\ell=\ell(X)$ is the length of the first irreducible part of $X$, and 
\item we have that 
\begin{equation}\label{E_renewal}
\frac{A_n^*}{nA_n}
=\E\left[
\frac{1}{\cI_n}
\right],
\end{equation}
where $\cI_n$ is the number of irreducible 
parts in a uniformly random structure
of length $n$. 
\end{enumerate}
\end{lemma}

It might be helpful to think of each $(X,m)$  
as encoding a unique structure $X(m)$ of length $n$, 
obtained by shifting $X$ by 
some magnitude $m$. 

\cref{L_renewal}
is proved in \cite{Kol23}, however, 
the following proof is 
simpler. In particular,  
\eqref{E_renewal} follows 
by the exchangeability of the irreducible parts. 

\begin{proof}
Let $B_n$ be the number of pairs $(X,m)$
as above. We will show that 
$nA_n=\sum_{k=1}^n B_k A_{n-k}$, 
as this implies $A_n^*=B_n$, 
as claimed in (1). 

Note that 
$nA_n$ counts $X$
of length $n$ with a marked point
$1\le j\le n$.  
If we split such an $X$ 
at the start of its 
irreducible part containing $j$ then, 
for some $1\le k \le n$, 
we obtain an unmarked structure of length $n-k$ 
and a structure of length $k$ 
with a mark in the first irreducible component. 
This procedure is injective, 
and its image is enumerated 
by $\sum_{k=1}^n B_k A_{n-k}$, 
proving the claim. 

Next, we will prove (2). 
By (1), we have that 
$A_n^*
=\sum_X\ell(X)$, summing over $X$
of length $n$. 
Hence $A_n^*/A_n$
is the expected length 
of the first irreducible component in 
a uniformly random $X$ of length $n$. 
Since the 
irreducible parts in such an $X$
are exchangeable, 
this equals 
$\E[n/\cI_n]$, 
and \eqref{E_renewal} follows. 
\end{proof}

As discussed in \cite{Kol23}, 
\eqref{E_renewal} gives 
probabilistic meaning to the right hand
side of \eqref{E_asySn}. Specifically, 
\[
\exp \left(-
\sum_{k=1}^\infty\frac{N_k}{k4^k}
\right)
=
\lim_{n\to\infty} \E\left[\frac{1}{\cI_n}\right]
\]
is the asymptotic expected inverse  
number of irreducible parts in a uniformly random
score sequence of length $n$.

\section{Strong score sequences}
We observe that 
$S_n$ is a renewal sequence. 
The irreducible parts of a score sequence 
are separated by 
the points $k$ for which 
$\sum_{i=1}^k s_i={k\choose2}$. 
Indeed, as discussed in \cite{Moo68}, 
score sequences 
with only one irreducible part
(such that $\sum_{i=1}^k s_i>{k\choose2}$,  
for all $0<k<n$) 
are called {\it strong,} 
since a tournament with a 
strong score sequence is strongly connected. 

By \cref{L_renewal}, 
to prove 
\eqref{E_CDFS} 
we need to show that 
$N_n$ enumerates pairs $(S,m)$, 
where $S$ is a score sequence of length $n$ and $0\le m <\ell $, where 
$\ell=\ell(S)$ is the length of the first irreducible 
part of $S$. In what follows, 
we will give a simple geometric 
explanation for this relationship.

\begin{figure}[h!]
\centering
\includegraphics[scale=1]{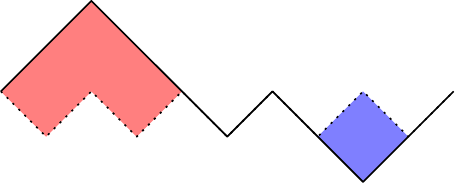}
\caption{A bridge $B$ (solid) of length $10=2\cdot 5$, 
with down steps at times 
$3$, $4$, $5$, $7$ and $8$. 
There are 3 diamonds (red) above and 1 diamond 
(blue) below
the sawtooth bridge (dotted), 
so its diamond area is
$a(B)=3-1=2$. 
Equivalently, in terms of its
down steps, 
$a(B)=-5^2+(3+4+5+7+8)=2$. 
}
\label{F_aB}
\end{figure}

\section{Diamond areas}
\label{S_areas}

We identify each sequence 
$1\le d_1<\cdots< d_n\le 2n$
with a bridge $B=(0=B_0,B_1,\ldots,B_{2n}=0)$
taking {\it down steps} $B_t-B_{t-1}=-1$
at times $t=d_i$, and {\it up steps}
$B_t-B_{t-1}=+1$ at all other times $1\le t\le 2n$. 

For reasons discussed below, 
our point of reference will be  the {\it sawtooth bridge}
$\check B=(0,-1,0,\ldots,-1,0)$,  
with down steps 
at odd times $d_i=2i-1$ 
and up steps at even times. 

For a bridge $B$ of length $2n$, 
we let $a(B)$ be  
$1/2$ of the area of $B$
above $\check B$, 
calculated as follows: 
\begin{align}\label{E_aB}
a(B)
&=\frac{1}{2}\sum_{t=0}^{2n}(B_t-\check B_t)
=\frac{1}{2}
\left[n+\sum_{t=1}^{2n}(2n+1-t)(B_t-B_{t-1})\right]\nonumber\\
&=-n^2+\sum_{i=1}^n d_i.
\end{align}
We call $a(B)$ the {\it diamond area}
of $B$. 
Graphically, 
$a(B)$ is the 
signed number of 
{\it diamonds} (rotated squares) 
between $B$ and $\check B$, 
as in \cref{F_aB}.

Crucially, 
we note that $a(B)\equiv \sum_{i=1}^n d_i$ mod $n$.

There are $2N_n$ bridges $B$ of length $2n$
with $a(B)\equiv0$ mod $n$. 
Indeed, 
such bridges that  furthermore end with an up step
correspond to sequences with 
$d_n\le 2n-1$ and are enumerated by $N_n$. 
Reflecting any such bridge $B$
over the $x$-axis yields a bridge $B'$
with $d_n=2n$ and 
$a(B')=n-a(B)$.

Following Erd\H{o}s and Moser (see 
\cite{Moo68}), we associate each score sequence
$0\le s_1\le\cdots\le s_n\le n-1$ 
with the bridge $B$ taking down steps
at times $d_i=s_i+i$. 
Informally, this bridge is obtained by drawing
the bar graph of the score sequence, and then 
rotating clockwise by $\pi/4$. 

Since $\sum_i s_i={n\choose2}$, 
it follows that $\sum_i d_i=n^2$, and 
so $a(B)=0$ for each such $B$. 
In fact, $B$ corresponds to a score sequence
if and only if 
$a(B)=0$ and  $a(B^{(2k)})\ge0$, for all 
sub-bridges  $B^{(2k)}=(B_0,B_1,\ldots,B_{2k})$ 
of $B$ with $B_{2k}=0$, 
since $a$ is monotone between such times.
This is simply a rephrasing of 
Landau's theorem  \cite{Lan53}
in terms of bridges. 

The sawtooth bridge $\check B$ is associated with  
score sequence $(0,1,\ldots,n-1)$, whose bar graph
is a ``staircase.'' The reason for the choice of $\check B$, 
in the definition of $a(B)$ in 
\eqref{E_aB}
above, is that $(0,1,\ldots,n-1)$ is extremal, in the sense that, by 
Landau's theorem, it has minimal partial sums ${k\choose 2}$. 
(Geometrically, 
$(0,1,\ldots,n-1)$ is a vertex of 
the polytope $\Pi_{n-1}$ discussed in 
\cref{E_perm}.)

\section{Kleitman's intuition}
\label{S_Kleitman}

Kleitman \cite{Mos71} (cf.\ \cite{Kle70,WK83}) 
observed that $S_n$ can be bounded by 
cyclically shifting the positive/negative areas 
enclosed by the sawtooth bridge $\check B$ 
and bridges 
$B$ with $a(B)=0$. 
By Raney \cite{Ran60}, 
this procedure can shift 
any such $B$
into a bridge $B'$
associated with a score sequence. 
The difficulty is that the shift is not unique. 
To bound $S_n$, 
Kleitman considered  
the average number of such
shifts in a random $B$.

\begin{figure}[h!]
\centering
\includegraphics[scale=0.95]{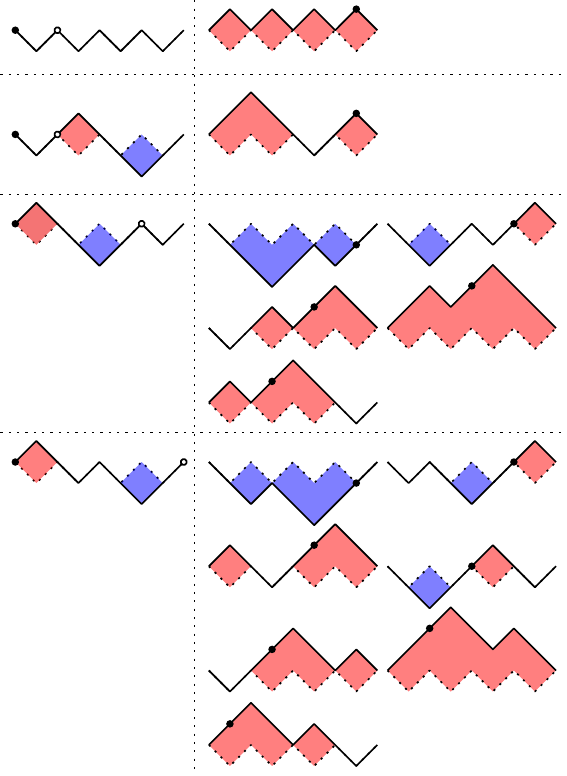}
\caption{The $S_4=4$ bridges $B$
(at left) 
associated with  the 
score sequences 0123, 0222, 1113 
and 1122
and the $2N_4=4+14=18$ bridges
$B'$ (at left and right) with 
$a(B)\equiv0$ mod $4$.
To obtain a bijective correspondence, 
we cyclically shift bridges $B$
by some $m$ (black dots) less than 
the length $2\ell$ (white dots) of their 
first irreducible parts. 
}
\label{F_n4}
\end{figure}

\section{Proof that $S_n^*=N_n$}
\label{S_proof}

As discussed, Kleitman studied bridges 
$B$ with $a(B)=0$ in order to bound
the asymptotics of $S_n$. 
In this section, we relate the precise asymptotics 
of $S_n$ to cyclic shifts of bridges $B$ with 
$a(B)\equiv 0$ mod $n$.

Specifically, 
using \cref{L_renewal} we show that 
$2S_n^*=2N_n$, by identifying a 
simple bijection $\phi$ that assigns each of the $2N_n$
many bridges 
$B$ with $a(B)\equiv 0$ mod $n$
to a unique shift of a $B'$
associated with a score sequence of length $n$. 

The bijection is quite natural,
as depicted in \cref{F_n4}.
The key observation is that
if we translate the sawtooth bridge 
$\check B$ 
by some $\delta\in\Z$, the diamond area $a$ 
of a bridge of length $2n$
becomes 
$a'=a-\delta n$;
see \eqref{E_aB}. 

Suppose that a bridge $B'$ 
corresponds to a score sequence, 
and that its first irreducible part is of length $\ell$.
Then, for $0\le m<2\ell$, we let 
$\phi(B',m)$ be the bridge $B$, 
with $a(B)\equiv 0$ mod $n$, 
obtained by
cyclically shifting the increments 
of $B$
to the left by $m$.

On the other hand, suppose that $a(B)\equiv 0$ mod $n$. 
The inverse bijection $\phi^{-1}(B)=(B',m)$
is obtained as follows: First, we find the unique 
shift of $\check B$ by some $\delta$ that makes the diamond
area equal to $a'=0$. Then, along this shifted
sawtooth path, we find the rightmost point, 
intersected by some vertical line $x=2n-m$, 
such that the bridge started from this point
has non-negative cumulative diamond areas
(with respect to the shifted $\check B$). 
Such a point 
exists by Raney \cite{Ran60}.

By \cref{L_renewal}(1) it follows that 
$2S_n^*=2N_n$, and hence 
$S_n^*=N_n$.

\section{Moser's conjecture}
\label{S_OurCase}

Combining \eqref{E_CDFS}, \eqref{E_Nn} and 
\eqref{E_asyAn}, we obtain the following short proof 
of Moser's conjecture. 
By \eqref{E_Nn}, $N_n\sim {2n\choose n}/2n$, so, 
in particular, by \eqref{E_CDFS} and Stirling's approximation, 
$(S_n/4^n)^*=N_n/4^n$ is regularly varying with index $\gamma=-3/2$. 
Therefore, \cref{T_main} follows by \eqref{E_CDFS}. 

On the other hand, past attempts 
\cite{Moo68,Kle70,Mos71,WK83,Tak86,Tak91,KP00}
at a  
direct analysis culminated with 
$S_n=\Theta(4^n/n^{5/2})$. 
However, let us mention that the recent work 
by the second and third authors \cite{DK24a}
shows that a direct 
approach (as outlined by Kleitman \cite{Kle70}) is possible. 
This approach  
yields additional insights, such as the Airy integral
\cite[Corollary 3]{DK24a}
and the scaling limit \cite[Theorem 5]{DK24a}, 
but is considerably more involved.


\makeatletter
\renewcommand\@biblabel[1]{#1.}
\makeatother

\providecommand{\bysame}{\leavevmode\hbox to3em{\hrulefill}\thinspace}
\providecommand{\MR}{\relax\ifhmode\unskip\space\fi MR }
\providecommand{\MRhref}[2]{%
  \href{http://www.ams.org/mathscinet-getitem?mr=#1}{#2}
}
\providecommand{\href}[2]{#2}

\end{document}